# ON THE SIGNATURES OF EVEN 4–MANIFOLDS

CHRISTIAN BOHR

ABSTRACT. In this paper, we prove a number of inequalities between the signature and the Betti numbers of a 4–manifold with even intersection form and prescribed fundamental group. Furthermore, we introduce a new geometric group invariant and discuss some of its properties.

## 1. Introduction

In this paper, we investigate which even forms can be realized as intersection forms of smooth 4–manifolds with a prescribed fundamental group. In particular, we prove estimates between the signature and the Betti numbers of a 4–manifold with even intersection form which are implied by special algebraic properties of the fundamental group.

If $X$ is a 4–manifold with non–trivial even intersection form $Q$, then $Q$ is indefinite (see [Do]). According to the Hasse–Minkowski classification (see for instance [MH]), it must therefore be of the form

$$Q = pE_8 \oplus qH,$$

where $E_8$ is the negative definite bilinear form defined by the Dynkin diagram of the Lie algebra $E_8$, $H$ is the hyperbolic form of rank 2, and $q > 0$. Several people have stated the following

**Conjecture:** If $X$ is a 4–manifold with even intersection form and signature $\tau(X)$, then $\frac{5}{4}|\tau(X)| \le b_2(X)$.

We will refer to this as the $\frac{5}{4}$–conjecture. Note that, with the notation above, the conjecture reads $|p| \le q$. Furuta has proved this conjecture in the spin case (see [Fu]). If $X$ is a 4–manifold with even intersection form and there is no 2–torsion in $H_1(X;\mathbb{Z}) = H_1(\pi_1(X);\mathbb{Z})$, then $X$ is spin and the conjecture is true by Furuta's Theorem. So the only interesting cases are 4–manifolds with 2–torsion in their first homology. Using covering arguments, we will prove that the $\frac{5}{4}$–conjecture is true for some special fundamental groups, including all finite and abelian groups (Corollary 2).

In Section 2, we show that every 4–manifold with even intersection form has a finite cover which is spin (Theorem 1), and determine the minimal degree of

1991 *Mathematics Subject Classification.* Primary 57M05 ; Secondary 57R19.

The author has been supported by the Graduiertenkolleg "Mathematik im Bereich ihrer Wechselwirkung mit der Physik" at the University of Munich.

such a covering (Theorem 2). In Section 3, we use this result to estimate the signature of a 4–manifold which has even intersection form or is finitely covered by a 4–manifold with even intersection form. In particular, we prove a $\mathbb{Z}_2$–version of Furuta's Theorem for these 4–manifolds (Theorem 3) and an $L^2$–version of the $\frac{5}{4}$–conjecture for residually finite fundamental groups (Theorem 5).

It turns out that, for special fundamental groups, for example for amenable groups, we have inequalities between the signature and the Euler characteristic which are in many cases considerably stronger than the $\frac{5}{4}$–conjecture (Theorem 4). In order to investigate these relations between the characteristic numbers and the fundamental groups of even 4–manifolds systematically, we introduce in Section 4 a new geometric group invariant $r$, in the spirit of the invariants $p$ and $q$ defined by Kotschick in [Ko1] and Hausmann–Weinberger in [HW]. We compute the values of this invariant for a number of examples where the values of $p$ and $q$ are also known. Furthermore, we prove the existence of an infinite family of groups for which all the three invariants are different (Theorem 6).

During the preparation of this manuscript, the author was informed that R. Lee and T.J. Li have independently obtained some of the results contained in this paper, namely Theorem 1 and Corollary 1. Moreover, they give a proof of the $\frac{5}{4}$–conjecture if the 2–torsion of $H_1(X;\mathbb{Z})$ equals $\mathbb{Z}_2 \oplus \mathbb{Z}_2$ (see [LL]).

It is clear that the results presented in this paper immediately yield estimates for the minimal genera of surfaces in 4–manifolds. Suppose for example that $F \subset X$ is an embedded surface in a 4–manifold such that the complement $X \setminus F$ has amenable fundamental group. Taking a branched cover yields a 4–manifold $Y$ whose signature and Euler characteristic can be computed from the genus of $F$ and its homology class. Under certain conditions at the homology class of $F$, the manifold $Y$ will have even intersection form or will even be spin. Furthermore, it is easy to see that the fundamental group of $Y$ is again amenable, and hence one can use the results of Section 3 to derive an estimate for the genus of $F$. Details can be found in [Bo].

Finally, I would like to thank my advisor, Dieter Kotschick, for many helpful suggestions and discussions.

## 2. Spin coverings

In this section, we will prove that every 4–manifold $X$ with even intersection form has a finite cover which is spin. It turns out that the degree of this covering can be chosen to be a power of two, more precisely the exponent of the 2–torsion of $H_1(X;\mathbb{Z})$. Furthermore, we will see that, in general, the covering we construct has minimal degree among all spin coverings. First we need some technical Lemmas.

**Lemma 1:** *Let $X$ be a 4–manifold with second Stiefel–Whitney class $w_2(X) \in H^2(X;\mathbb{Z}_2)$ and consider the inclusion $Ext(H_1(X);\mathbb{Z}_2) \hookrightarrow H^2(X;\mathbb{Z}_2)$ given by the Universal Coefficient Theorem. Then the following conditions are equivalent:*

1. *$X$ has even intersection form*

2. $w_2(X) \in Ext(H_1(X); \mathbb{Z}_2)$.

*Proof.* Consider the diagram

$$\begin{array}{ccccc} Ext(H_1(X;\mathbb{Z});\mathbb{Z}) & \longrightarrow & H^2(X;\mathbb{Z}) & \xrightarrow{\kappa_1} & Hom(H_2(X;\mathbb{Z});\mathbb{Z}) \\ \rho_1 \downarrow & & \rho_2 \downarrow & & \rho_3 \downarrow \\ Ext(H_1(X;\mathbb{Z});\mathbb{Z}_2) & \longrightarrow & H^2(X;\mathbb{Z}_2) & \xrightarrow{\kappa_2} & Hom(H_2(X;\mathbb{Z});\mathbb{Z}_2) \end{array}$$

where the $\kappa_i$ are given by evaluation and the $\rho_i$ by reduction modulo 2. The Universal Coefficient Theorem states that the rows are exact. Choose a characteristic class $c \in H^2(X;\mathbb{Z})$. Now suppose that $X$ has even intersection form. Then

$$\kappa_1(c)(x) = Q_X(c, PD(x)) \equiv Q_X(PD(x), PD(x)) \equiv 0 \mod 2$$

for every $x \in H_2(X;\mathbb{Z})$, therefore $0 = \rho_3\kappa_1(c) = \kappa_2(w_2(X))$, and the exactness implies $w_2(X) \in Ext(H_1(X);\mathbb{Z}_2)$. This proves $1 \implies 2$. For the proof of $2 \implies 1$ suppose that $w_2(X) \in Ext(H_1(X);\mathbb{Z}_2)$, i.e. $\kappa_2(w_2(X)) = 0 = \rho_3\kappa_1(c)$. This implies

$$Q_X(c, x) \equiv \kappa_1(c)(x) \equiv 0 \mod 2$$

for every $x \in H^2(X;\mathbb{Z})$, and therefore $X$ has even intersection form. □

**Definition 1:** *A covering map $p : \tilde{X} \to X$ with connected $\tilde{X}$ will be called* **abelian** *if it is regular (i.e. the group of deck transformations acts transitively on each fibre or, equivalently, the characteristic subgroup $\pi_1(\tilde{X}) \subset \pi_1(X)$ is a normal subgroup) and the group of deck transformations is abelian. It will be called* **cyclic** *if it is abelian and the group of deck transformations is cyclic.*

**Lemma 2:** *Suppose that $G$ is a finite abelian group with 2–primary part $G(2)$ and $2^\mu$ the exponent of $G(2)$. Then there is for every $w \in Ext(G;\mathbb{Z}_2)$ a homomorphism $f : G \to \mathbb{Z}_{2^\mu}$ such that the induced map $Ext(f,\mathbb{Z}_2)$ maps the non–zero element of $Ext(\mathbb{Z}_{2^\mu},\mathbb{Z}_2) = \mathbb{Z}_2$ to $w$.*

*Proof.* There is a decomposition $G = G(2) \oplus H$, where the order of $H$ is odd, and from the exact sequence

$$0 \longrightarrow G(2) \longrightarrow G \longrightarrow H \longrightarrow 0$$

we can deduce the exact sequence

$$Ext(H;\mathbb{Z}_2) \longrightarrow Ext(G;\mathbb{Z}_2) \longrightarrow Ext(G(2);\mathbb{Z}_2) \longrightarrow 0$$

(see [Mc], Chapter III, 3.2 and 3.7). But $Ext(H;\mathbb{Z}_2) = 0$ and hence the inclusion induces an isomorphism $Ext(G;\mathbb{Z}_2) = Ext(G(2);\mathbb{Z}_2)$. Therefore we can assume that $G = G(2)$ is a 2–group. There is a unique decomposition $G = \bigoplus_{i=1}^r \mathbb{Z}_{2^{\mu_i}}$ of $G$ as a sum of cyclic groups, and $\mu = \max_i\{\mu_i\}$. A homomorphism $f : G \to \mathbb{Z}_{2^\mu}$ is then determined by the images of the generators and therefore can be described by a vector $(a_1, \ldots, a_r)$ with $0 \leq a_i < 2^\mu$ subject to the conditions $2^\mu | 2^{\mu_i} a_i$ for all $i$. Now assume a homomorphism $f$ is given and consider the standard resolution

$$0 \longrightarrow \mathbb{Z}^r \xrightarrow{\iota} \mathbb{Z}^r \xrightarrow{\pi} G \longrightarrow 0$$

where $\iota(x_1, \ldots, x_r) = (2^{\mu_1} x_1, \ldots, 2^{\mu_r} x_r)$. Define homomorphisms $f_i : \mathbb{Z}^r \to \mathbb{Z}$, $i = 1, 2$, by

$$f_1(x_1, \ldots, x_r) = \sum_i x_i a_i,$$

$$f_2(x_1, \ldots, x_r) = \sum_i 2^{\mu_i - \mu} a_i x_i.$$

Note that the last sum makes sense, since $a_i 2^{\mu_i} = 0 \in \mathbb{Z}_{2^\mu}$, and therefore $a_i 2^{\mu_i}$ is divisible by $2^\mu$ in $\mathbb{Z}$. Then the diagram

$$\begin{array}{ccccccccc}
0 & \longrightarrow & \mathbb{Z}^r & \xrightarrow{\iota} & \mathbb{Z}^r & \xrightarrow{\pi} & G & \longrightarrow & 0 \\
& & \downarrow {\scriptstyle f_2} & & \downarrow {\scriptstyle f_1} & & \downarrow {\scriptstyle f} & & \downarrow \\
0 & \longrightarrow & \mathbb{Z} & \xrightarrow{j} & \mathbb{Z} & \longrightarrow & \mathbb{Z}_{2^\mu} & \longrightarrow & 0
\end{array}$$

commutes, where $j(x) = 2^\mu x$. Now, using the definition of $Ext(\cdot; \mathbb{Z}_2)$ and the fact that the induced maps $\iota^* : Hom(\mathbb{Z}^r, \mathbb{Z}_2) \to Hom(\mathbb{Z}^r, \mathbb{Z}_2)$ and $j^* : Hom(\mathbb{Z}, \mathbb{Z}_2) \to Hom(\mathbb{Z}, \mathbb{Z}_2)$ are both zero, we obtain that there are isomorphisms $Ext(G; \mathbb{Z}_2) = Hom(\mathbb{Z}^r, \mathbb{Z}_2)$ and $Ext(\mathbb{Z}_{2^\mu}, \mathbb{Z}_2) = Hom(\mathbb{Z}; \mathbb{Z}_2)$ fitting into the diagram

$$\begin{array}{ccc}
Hom(Z^r, \mathbb{Z}_2) & \xrightarrow{\cong} & Ext(G; \mathbb{Z}_2) \\
{\scriptstyle f_2^*} \uparrow & & \uparrow {\scriptstyle Ext(f; \mathbb{Z}_2)} \\
Hom(\mathbb{Z}, \mathbb{Z}_2) & \xrightarrow{\cong} & Ext(\mathbb{Z}_{2^\mu}, \mathbb{Z}_2)
\end{array}$$

This gives us a description of $Ext(f; \mathbb{Z}_2)$ as a map from $Hom(\mathbb{Z}, \mathbb{Z}_2) = \mathbb{Z}_2$ to $Hom(\mathbb{Z}^r, \mathbb{Z}_2) = \mathbb{Z}_2^r$, and we see that the image of the nonzero element is given by the vector $(2^{\mu_i - \mu} a_1, \ldots, 2^{\mu_i - \mu} a_r) \in \mathbb{Z}_2^r$ (here again we use that $a_i 2^{\mu_i}$ is divisible by $2^\mu$). Now $w \in Ext(G, \mathbb{Z}_2)$ is described by a vector $(w_1, \ldots, w_r)$, and if we choose $f$ to be the homomorphism defined by the vector $(2^{\mu - \mu_i} w_1, \ldots, 2^{\mu - \mu_i} w_r)$, then $Ext(f; \mathbb{Z}_2)(1) = w$. $\square$

**Lemma 3:** *Let $X$ be a 4–manifold with $H_1(X; \mathbb{Z}) = \mathbb{Z}_{2^k}$, $k \geq 2$, and suppose that $\tilde{X} \to X$ is a connected double cover. Then $X$ is spin if and only if $\tilde{X}$ is spin.*

*Proof.* The proof is a slight modification of an argument that appeared in a paper by Massey (see [Ma]). We have to show that the map $H^2(X; \mathbb{Z}_2) \to H^2(\tilde{X}; \mathbb{Z}_2)$ induced by the projection is one–to–one. For this purpose, consider the exact sequence

$$H^1(X; \mathbb{Z}_2) \xrightarrow{\cup w_1} H^2(X; \mathbb{Z}_2) \longrightarrow H^2(\tilde{X}; \mathbb{Z}_2)$$

which is a part of the Gysin sequence associated to $\tilde{X} \to X$. Since the covering is assumed to be non-trivial, $w_1 \neq 0$, and since $H^1(X; \mathbb{Z}_2) = \mathbb{Z}_2$, it is sufficient to prove $w_1 \cup w_1 = 0$. But $w_1 \cup w_1 = \beta(w_1)$, where $\beta$ is the boundary in the long exact sequence associated to the short exact sequence $\mathbb{Z}_2 \to \mathbb{Z}_4 \to \mathbb{Z}_2$ of coefficient groups. Therefore we only have to check that the homomorphism $H^1(X; \mathbb{Z}_4) \to H^1(X; \mathbb{Z}_2)$ induced by the projection $\mathbb{Z}_4 \to \mathbb{Z}_2$, is onto, or, in other words, that the canonical map $Hom(\mathbb{Z}_{2^k}; \mathbb{Z}_4) \to Hom(\mathbb{Z}_{2^k}; \mathbb{Z}_2)$ is onto, which is clearly the case if and only if $k \geq 2$. □

The following Theorem is the main technical result in this section and asserts in particular that every 4–manifold with even intersection form has a finite spin covering.

**Theorem 1:** *If $X$ is a connected closed and oriented 4–manifold with even intersection form, then there is a finite cyclic cover $\tilde{X} \to X$ which is spin. Furthermore, the degree of this cover is a power of $2$ and at most the exponent of the 2–primary part of $\operatorname{Tor} H_1(X; \mathbb{Z})$.*

*Proof.* Assume that $X$ has even intersection form, but is not spin, i.e. $w_2(X) \neq 0$. Using $Ext(F; \mathbb{Z}_2) = 0$ for a finitely generated free abelian group $F$ and the exact sequence (see [Mc], Chapter III, 3.2 and 3.7) induced by a splitting $H_1(X; \mathbb{Z}) = F \oplus T$ into a free part and the torsion, we see that the inclusion $\operatorname{Tor} H_1(X; \mathbb{Z}) \to H_1(X; \mathbb{Z})$ induces an isomorphism $Ext(\operatorname{Tor} H_1(X); \mathbb{Z}_2) = Ext(H_1(X); \mathbb{Z}_2)$. This means that for every homomorphism $\varphi : H_1(X; \mathbb{Z}) \to G$ into some group $G$, $Ext(\varphi; \mathbb{Z}_2)$ is determined by the restriction of $\varphi$ to the torsion subgroup of $H_1(X; \mathbb{Z})$. Now let $2^\mu$ denote the exponent of the 2–primary part of $\operatorname{Tor} H_1(X; \mathbb{Z})$. By Lemma 2, there is a homomorphism $\varphi : H_1(X; \mathbb{Z}) \to \mathbb{Z}_{2^\mu}$ such that $Ext(\varphi; \mathbb{Z}_2)$ maps the non–zero element in $Ext(\mathbb{Z}_{2^\mu}; \mathbb{Z}_2)$ to $w_2(X) \in Ext(H_1(X; \mathbb{Z}); \mathbb{Z}_2)$ (if $b_1(X) > 0$, we can choose this homomorphism to be zero on $F$). In general, this homomorphism will not be surjective, but its image is a cyclic subgroup $\mathbb{Z}_{2^\nu} \subset \mathbb{Z}_{2^\mu}$. Note that, since $w_2(X) \neq 0$, $\varphi \neq 0$, and therefore $\nu \geq 1$. A simple calculation shows that the inclusion induces an isomorphism $Ext(\mathbb{Z}_{2^\mu}; \mathbb{Z}_2) = Ext(\mathbb{Z}_{2^\nu}; \mathbb{Z}_2)$, hence we can assume that $\mu = \nu$ and $\varphi$ is a surjection. Let $\tilde{X} \to X$ denote the connected cover defined by $\varphi$ and consider the exact sequence

$$H_1(\tilde{X}; \mathbb{Z}) \xrightarrow{\pi_*} H_1(X; \mathbb{Z}) \xrightarrow{\varphi} \mathbb{Z}_{2^\mu} \longrightarrow 0.$$

Applying the Ext–functor shows that the composition

$$Ext(\mathbb{Z}_{2^\mu}; \mathbb{Z}_2) \xrightarrow{\varphi^*} Ext(H_1(X; \mathbb{Z}); \mathbb{Z}_2) \longrightarrow Ext(A; \mathbb{Z}_2)$$

is zero, and since $w_2(X)$ is lying in the image of $\varphi^*$, this implies that $Ext(\pi_*; \mathbb{Z}_2)$ maps $w_2(X)$ to zero, and therefore $\tilde{X}$ is spin. □

It is of course tempting to conjecture that there is always a double cover which is spin. However, the next result asserts that the covering constructed in the proof of Theorem 1 has in general the lowest possible degree.

**Theorem 2:** *For every finitely generated abelian group $G$, there is a 4–manifold $X$ with even intersection form and $H_1(X;\mathbb{Z}) = G$ such that for no covering $\tilde{X} \to X$ of degree less than the exponent of the 2–torsion of $G$, the covering space $\tilde{X}$ is spin.*

*Proof.* First assume that $G = \mathbb{Z}_n$ where $n = 2^\mu$ is a power of 2. Pick a rational homology 4–sphere $\Sigma_n$ with fundamental group $\mathbb{Z}_n$ which is not spin (such a manifold exists by Proposition 1 in [HK], and there is also an elementary construction using surgery on $S^1 \times S^3$). Now assume that $X \to \Sigma_n$ is a finite cover of degree $m = 2^\nu$ less than $n$. Then this covering corresponds to a subgroup $\mathbb{Z}_{2^{\mu-\nu}} \subset \mathbb{Z}_{2^\mu}$ and hence can be decomposed into a sequence

$$X = X_{2^{\mu-\nu}} \longrightarrow X_{2^{\mu-\nu+1}} \longrightarrow \cdots \longrightarrow X_{2^\mu} = \Sigma_n$$

of double coverings, corresponding to the sequence $\mathbb{Z}_{2^{\mu-\nu}} \subset \cdots \subset \mathbb{Z}_{2^\mu}$. By Lemma 3, applied inductively to these coverings, we obtain that $Y$ is not spin if $m < n$. This proves the assertion in the special case $G = \mathbb{Z}_{2^\mu}$. Now suppose we are given any finitely generated abelian group $G$ and let $e = 2^\mu$ denote the exponent of the 2–torsion $\text{Tor}_2(G)$. We can decompose the 2–torsion into a sum of cyclic 2–groups, and, by definition of the exponent, one of the summands is $\mathbb{Z}_{2^\mu}$. Hence we have a decomposition $G = H \oplus \mathbb{Z}_{2^\mu}$. Now pick an arbitrary spin 4–manifold $Y$ with fundamental group $H$ (which one can construct using 1–surgery on a connected sum of copies of $S^1 \times S^3$) and let $X = Y \# \Sigma_n$. If we have a covering $\tilde{X} \to X$ of degree less than $e$, the covering space splits into a fibre sum of a covering of $Y$ and a covering of $\Sigma_n$, both possibly not connected. If $\tilde{X}$ is spin, each connected component of the covering over $\Sigma_n$ is also spin. But every component is itself a covering space over $\Sigma_n$ of degree less than $e$ and, by the arguments used in the special case $G = \mathbb{Z}_{2^\mu}$, is therefore not spin, a contradiction. □

**Remark 1:** A part of the assertion of Theorem 1, namely the fact that an even 4–manifold is finitely covered by a spin manifold, can also be proved using a result of Deligne and Sullivan. Suppose $X$ has even intersection form. From the proof of Lemma 1, one can deduce that there is a characteristic class $c \in H^2(X;\mathbb{Z})$ which is torsion. In other words, there is a $\text{Spin}^{\mathbb{C}}$-structure on $X$ such that a certain tensor power of its associated line bundle is trivial, and its real Chern class is zero. This implies that there is a flat connection on this bundle, hence this bundle is associated to the universal covering of $X$ via the monodromy $\pi_1(X) \to U(1)$. Deligne and Sullivan prove in [DS] that under these conditions, there is a finite cover on which this bundle is trivial. However, this proof gives no detailed information on the covering, at least not at first glance.

Conversely, a special case of the Theorem of Deligne and Sullivan can also be proved using arguments similar to those given in the proof of Theorem 1. Suppose that $X$ is a compact smooth manifold and $E \to X$ a flat complex line bundle, in the sense that there is a flat U(1)–connection on $E \to X$. Then the real Chern class $c_1(E)_{\mathbb{R}}$ is zero, and therefore the Chern class $c_1(E) \in H^2(X;\mathbb{Z})$ must be torsion. By the Universal Coefficient Theorem, the torsion of $H^2(X;\mathbb{Z})$ is $Ext(H_1(X;\mathbb{Z});\mathbb{Z})$. Now choose a splitting $H_1(X;\mathbb{Z}) = F \oplus T$ into a free abelian group $F$ and a finite group $T$ and consider the covering $\tilde{X} \to X$ given by projection onto $T$. Then the quotient $A = H_1(\tilde{X};\mathbb{Z})/ker\pi_*$ is isomorphic to $F$, hence free, and $Ext(A;\mathbb{Z}) = 0$. But the map $Ext(H_1(X;\mathbb{Z});\mathbb{Z}) \to Ext(H_1(\tilde{X};\mathbb{Z});\mathbb{Z})$ factors through $Ext(A;\mathbb{Z})$ and is therefore zero. This means that the whole torsion of $H^2(X;\mathbb{Z})$ becomes zero under the pullback to $H^2(\tilde{X};\mathbb{Z})$ and the pullback of $E$ is trivial.

In the case of a flat complex vector bundle of higher rank, this argument only proves that there is a finite cover on which all Chern classes are zero, but this does not imply that the bundle is trivial (a simple example is the non–trivial $\mathbb{C}^2$–bundle over $S^5$ given by the non–zero element in $\pi_4(U(2)) = \mathbb{Z}_2$, this bundle is not trivial although all its Chern classes are zero).

## 3. Signature estimates for even 4–manifolds

In this section, we use covering arguments to obtain a number of inequalities between the signature and the Betti numbers of a 4–manifold having even intersection form. As a first result, we obtain a $\mathbb{Z}_2$–version of Furuta's Theorem for 4–manifolds with even intersection form. Furthermore, we will see that there are several classes of fundamental groups for which one can obtain non–trivial inequalities between the signature and the Euler characteristic of an even 4–manifold. But first, we need an estimate for the Betti numbers of a finite covering.

**Lemma 4:** *Let $X$ be a connected finite CW–complex and $\pi : \tilde{X} \to X$ be a connected abelian cover of degree $m = 2^{\mu}$. Then*

$$dim_{\mathbb{Z}_2} H^1(\tilde{X};\mathbb{Z}_2) \leq m\, dim_{\mathbb{Z}_2} H^1(X;\mathbb{Z}_2) - m + 1.$$

*Proof.* By induction on $\mu$. In the case $\mu = 0$, there is nothing to prove. Now assume that the Lemma is proved for $\mu - 1$ and $\tilde{X} \to X$ has degree $m = 2^{\mu}$. Choose a covering transformation $\sigma$ of order 2 and let $Y = \tilde{X}/\sigma$. We obtain a double cover $\tilde{X} \to Y$ and an abelian cover $Y \to X$ of degree $2^{\mu-1}$. Now we have the exact sequence (Gysin sequence, see [Sp], Chapter 5, Section 7, Thm.11 )

$$H^0(Y;\mathbb{Z}_2) \xrightarrow{\cup w_1} H^1(Y;\mathbb{Z}_2) \longrightarrow H^1(\tilde{X};\mathbb{Z}_2) \longrightarrow H^1(Y;\mathbb{Z}_2)$$

from which we can deduce $dim_{\mathbb{Z}_2} H^1(\tilde{X};\mathbb{Z}_2) \leq 2\, dim_{\mathbb{Z}_2} H^1(Y;\mathbb{Z}_2) - 1$ (note that, since the covering is non–trivial, $w_1 \neq 0$). By the induction hypothesis, we have

the inequality

$$dim_{\mathbb{Z}_2} H^1(Y; \mathbb{Z}_2) \leq \frac{m}{2} dim_{\mathbb{Z}_2} H^1(X; \mathbb{Z}_2) - \frac{m}{2} + 1.$$

Substituting this into the estimate from the Gysin sequence, we obtain

$$dim_{\mathbb{Z}_2} H^1(\tilde{X}; \mathbb{Z}_2) \leq 2 \left( \frac{m}{2} dim_{\mathbb{Z}_2} H^1(X; \mathbb{Z}_2) - \frac{m}{2} + 1 \right) - 1$$
$$= m \, dim_{\mathbb{Z}_2} H^1(X; \mathbb{Z}_2) - m + 1,$$

as claimed. □

**Lemma 5:** *Let $X$ be a 4–manifold and let $t = dim_{\mathbb{Z}_2}(\text{Tor}_2 H_1(X; \mathbb{Z}) \otimes \mathbb{Z}_2)$ denote the number of summands in the (unique) decomposition of the 2–primary part $\text{Tor}_2 H_1(X; \mathbb{Z})$ of the 2–torsion of $H_1(X; \mathbb{Z})$ into cyclic subgroups. Then we have*

$$dim_{\mathbb{Z}_2} H^2(X; \mathbb{Z}_2) = b_2(X) + 2\,t.$$

*In particular, the number $dim_{\mathbb{Z}_2} H^2(X; \mathbb{Z}_2)$ is even if $X$ has even intersection form and the difference $dim_{\mathbb{Z}_2} H^2(X; \mathbb{Z}_2) - b_2(X)$ depends only on the fundamental group.*

*Proof.* This is an immediate consequence of Poincaré duality and the Universal Coefficient Theorem. □

**Lemma 6:** *Let $\pi : \tilde{X} \to X$ be a covering map of odd degree, where $\tilde{X}$ and $X$ are closed connected and oriented 4–manifolds.*
1. *If $\tilde{X}$ has even intersection form, then the same is true for $X$.*
2. *If $\tilde{X}$ is spin, then $X$ is spin.*

*Proof.* Let us start with the first assertion: Assume $X$ is not even, i.e. there is a class $c \in H^2(X; \mathbb{Z})$ such that $Q_X(c, c) \equiv 1 \mod 2$. Let $m$ denote the degree of the covering. Then

$$Q_{\tilde{X}}(\pi^* c, \pi^* c) = m Q_X(c, c) \equiv 1 \mod 2,$$

in contradiction to the assumption that $\tilde{X}$ has even intersection form. This proves that $X$ must be even.

Now assume that $\tilde{X}$ is spin. Since, by Wu's Theorem, $w_2(X)$ is completely determined by the condition $w_2 \cup x = x \cup x$ for all $x \in H^2(X; \mathbb{Z}_2)$, the manifold $X$ is spin if and only if $x \cup x = 0$ for all $x \in H^2(X; \mathbb{Z}_2)$. Now assume that $X$ is not spin. Then there is a class $x \in H^2(X; \mathbb{Z}_2)$ with $x \cup x \neq 0$. But again this implies $(p^* x) \cup (p^* x) \neq 0$, since the mapping degree of the projection is odd (note that every class in $H^4(\cdot; \mathbb{Z}_2)$ has a lift to $H^4(\cdot, \mathbb{Z})$, and therefore $H^4(X; \mathbb{Z}_2) \to H^4(\tilde{X}; \mathbb{Z}_2)$ is an isomorphism), and this is a contradiction, since $\tilde{X}$ is spin. □

Now we are ready to prove a $\mathbb{Z}_2$–version of Furuta's Theorem which also holds for non–spin 4–manifolds with even intersection form.

**Theorem 3:** *Let $X$ be a 4–manifold with non–zero signature $\tau(X)$.*

1. *If the intersection form of $X$ is even, then*
$$\frac{5}{4}|\tau(X)| \leq \dim_{\mathbb{Z}_2} H^2(X;\mathbb{Z}_2) - 2.$$

2. *If there is a finite abelian cover $\tilde{X} \to X$ such that the intersection form of $\tilde{X}$ is even, then*
$$\frac{5}{4}|\tau(X)| < \dim_{\mathbb{Z}_2} H^2(X;\mathbb{Z}_2).$$

*Proof.* We will use the abbreviation $b_k(X;\mathbb{Z}_2) = \dim_{\mathbb{Z}_2} H^k(X;\mathbb{Z}_2)$. Let us prove assertion 1 first. Assume that $X$ has even intersection form. By Theorem 1, there is a finite cyclic cover $\tilde{X} \to X$ with degree $m = 2^\mu$ for some $\mu$ such that $\tilde{X}$ is spin. Since the signature is multiplicative under finite coverings, the signature $\tau(\tilde{X}) = m\,\tau(X)$ is not zero. Furthermore, the multiplicativity of the Euler characteristic and Lemma 4 yield the estimate

$$\begin{aligned}
b_2(\tilde{X};\mathbb{Z}_2) &= m\,e(X) - 2 + 2b_1(\tilde{X};\mathbb{Z}_2) \\
&\leq m\,(2 + b_2(X;\mathbb{Z}_2) - 2b_1(X;\mathbb{Z}_2)) + 2\,(m\,b_1(X;\mathbb{Z}_2) - m) \\
&= 2m + m\,b_2(X;\mathbb{Z}_2) - 2m = m\,b_2(X;\mathbb{Z}_2).
\end{aligned}$$

In particular we obtain $b_2(\tilde{X}) \leq m\,b_2(X;\mathbb{Z}_2)$. Hence Furuta's Theorem, applied to $\tilde{X}$, yields $\frac{5}{4}|\tau(X)| \leq b_2(X;\mathbb{Z}_2) - \frac{2}{m}$. Since the two numbers $\frac{5}{4}|\tau(X)|$ and $b_2(X;\mathbb{Z}_2)$ are both even (recall Lemma 5), the claimed estimate follows.

Now let us turn to the proof of statement 2. By assumption, there is a finite abelian cover $\tilde{X} \to X$ such that $\tilde{X}$ has even intersection form. Let $G$ denote the group of covering transformations and $m = \#G$ its order. If $m$ is a prime power, then $m$ is a power of two or $X$ itself is even (according to Lemma 6). If $m$ has more than one prime factor, write

$$G = G(2) \oplus \bigoplus_{p \neq 2} G(p)$$

where $G(p)$ denotes the p-primary part of $G$. Let $H = \bigoplus_{p \neq 2} G(p)$ and $Y = \tilde{X}/H$. Then the order of $H$ is odd, hence $Y$ has even intersection form, and we have a finite abelian cover $Y = \tilde{X}/H \to X = \tilde{X}/G$ with group $G/H = G(2)$. In both cases we have constructed an abelian cover with even intersection form such that the degree of the covering is a power of 2, hence we can assume that the degree of $\tilde{X} \to X$ itself is a power of 2. Let $m$ denote this degree. Now we can apply the same calculation as above to conclude $b_2(\tilde{X};\mathbb{Z}_2) \leq m\,b_2(X;\mathbb{Z}_2)$. Furthermore the signature of $\tilde{X}$ is $m\,\tau(X) \neq 0$. Applying assertion 1. to $\tilde{X}$ leads to

$$\frac{5}{4}m\,|\tau(X)| < b_2(\tilde{X};\mathbb{Z}_2) \leq m\,b_2(X;\mathbb{Z}_2),$$

and dividing by $m$ leads to the desired estimate. □

As an immediate consequence, we obtain the $\frac{5}{4}$–conjecture in a special case, which is also a part of the main result of [LL].

**Corollary 1:** *The $\frac{5}{4}$–conjecture is true for even 4–manifolds $X$ for which the 2–torsion of $H_1(X; \mathbb{Z})$ is cyclic.*

*Proof.* If the order of the torsion is odd, our manifold $X$ will be spin, and the claim follows from Furuta's Theorem. If the order is even, then the number $t$ introduced in the statement of Lemma 5 is 1 and hence $b_2(X) = dim_{\mathbb{Z}_2} H^2(X; \mathbb{Z}_2) - 2$. Therefore Theorem 3 implies the inequality $\frac{5}{4}|\tau(X)| \leq b_2(X)$. □

The proof of Theorem 3 shows that the difficult part is to control the growth of the first Betti number of covering spaces. Therefore it seems unlikely to find a proof of the $\frac{5}{4}$–conjecture which is based on coverings. However, there are a number of special cases where taking coverings yields inequalities between the signature and the Euler characteristic of a 4–manifold with even intersection form which are even stronger than the $\frac{5}{4}$–conjecture. A first family of groups where this happens is the class of amenable groups.

**Definition 2:** *A group $G$ is called amenable if there is a left invariant mean on the space of bounded functions on $G$, i.e. a linear map $m : L^\infty(G) \to \mathbb{C}$ with $m(1) = 1$, $m(f) \geq 0$ whenever $f \geq 0$ and $m(gf) = m(f)$ for every $f \in L^\infty(G)$, $g \in G$.*

All finite and solvable groups are amenable, subgroups and homomorphic images of amenable groups are amenable, and the class of amenable groups is closed under extensions (see [Gl]). The standard example for a group which is not amenable is the free group with two generators. We now have the following

**Theorem 4:** *Let $X$ be a 4–manifold with even intersection form. Denote by $\tau(X)$ the signature of $X$ and by $e(X)$ its Euler characteristic. If $X$ has amenable fundamental group, then*

$$\frac{5}{4}|\tau(X)| \leq e(X).$$

*In particular, the $\frac{5}{4}$–conjecture holds for $X$ as soon as $b_1(X) > 0$.*

*Proof.* The proof is based on a modification of Gromov's argument for the fact that a 4–manifold $X$ with amenable fundamental group fulfills the inequality $|\tau(X)| \leq e(X)$ (see [Gr1], p. 83).

First suppose that $X$ is spin. If $b_1(X) = 0$, then $e(X) = 2 + b_2(X) > b_2(X)$, and the estimate follows from Furuta's Theorem. Therefore we can suppose that $b_1(X) \geq 1$. Then one has a surjective map $\pi_1(X) \to \mathbb{Z}$ and, by composition, a surjective homomorphism $\pi_1(X) \to \mathbb{Z}_k$ for every $k$. These maps determine coverings $\tilde{X} \to X$ (with group $\mathbb{Z}$) and $X_k \to X$ (with group $\mathbb{Z}_k$). Now assume that $\frac{5}{4}|\tau(X)| \geq e(X)+1$. Since $\frac{5}{4}\tau(X)$ and $e(X)$ are even, one even has $\frac{5}{4}|\tau(X)| \geq e(X) + 2$. Euler characteristic and signature are multiplicative under coverings,

and therefore
$$\frac{5}{4}|\tau(X_k)| \geq e(X_k) + 2k.$$

But
$$\frac{5}{4}|\tau(X_k)| - e(X_k) = \frac{5}{4}|\tau(X_k)| - b_2(X_k) - 2 + 2b_1(X_k) < 2b_1(X_k),$$

since by Furuta's Theorem, we know that $\frac{5}{4}|\tau(X_k)| - b_2(X_k) \leq 0$. Therefore we obtain $b_1(X_k) \geq k$. Now we can follow Gromov's arguments to conclude that there is a surjection from a subgroup of $\pi_1(X)$ of finite index (namely $\pi_1(X_k)$ for large $k$) onto a group of the form $\mathbb{Z}_p * \mathbb{Z}_p$ for some prime $p \geq 3$, and this in turn implies that $\pi_1(X)$ is not amenable (since $\mathbb{Z}_p * \mathbb{Z}_p$ is not amenable for $p \neq 2$), which is a contradiction.

If $X$ is not spin, we can choose a finite cover $\tilde{X} \to X$ which is spin. Since Euler characteristic and signature are multiplicative under finite coverings, the inequality for this finite cover implies the inequality for $X$. □

**Remark 2:** In fact we proved a stronger assertion: if the fundamental group of $X$ is not "large", i.e. there is no homomorphism from a subgroup of finite index onto the free group with two generators, then $\frac{5}{4}|\tau(X)| \leq e(X)$ (note that, for $p \neq 2$, $\mathbb{Z}_p * \mathbb{Z}_p$ contains a free group of rank $\geq 2$ as a subgroup of finite index). Furthermore, since Euler characteristic and signature are multiplicative under finite covers, the estimate is still true if $X$ itself does not have even intersection form, but a finite cover is even.

**Remark 3:** Note that if the $\frac{11}{8}$–conjecture is true, then the same argument yields the inequality $\frac{11}{8}|\tau(X)| \leq e(X)$ for all 4–manifolds with even intersection form and amenable fundamental group.

**Corollary 2:** *The $\frac{5}{4}$–conjecture is true for manifolds with finite or abelian fundamental group.*

*Proof.* Suppose that $X$ is a 4–manifold with even intersection form and finite fundamental group. Then the universal covering $\tilde{X}$ is spin (this follows for example from Lemma 1 and the naturality of $Ext$ and $w_2$), and therefore $\frac{5}{4}|\tau(\tilde{X})| \leq b_2(\tilde{X})$ by Furuta's Theorem. Let $m = \#\pi_1(X)$ denote the degree of the universal covering. The multiplicativity of the Euler characteristic then implies that
$$b_2(\tilde{X}) = m\, e(X) - 2$$
$$= 2m + mb_2(X) - 2$$
$$< m\,(b_2(X) + 2).$$

Applying Furuta's Theorem to $\tilde{X}$ and dividing by $m$ therefore yields $\frac{5}{4}|\tau(X)| < b_2(X) + 2$. Since $b_2(X)$ and $\frac{5}{4}$ are both even numbers, the claimed inequality follows.

Now assume that the fundamental group is abelian. If the first Betti number $b_1(X) = 0$, then the fundamental group is finite and we are done. If $b_1(X) > 0$, then the claim follows from Theorem 4. □

We have seen (Theorem 3) that there is a $\mathbb{Z}_2$–version of the $\frac{5}{4}$–conjecture which holds for arbitrary fundamental group. It it interesting that – at least in some cases – there is also an $L^2$–version, which is a consequence of a result of Lück concerning the growth of Betti numbers of certain covering spaces (see [Lu1]). We refer the reader who is not familiar with $L^2$–cohomology to [Lu2] for a nice introduction into this topic. Recall that a group $G$ is called **residually finite** if for every $g \in G \setminus \{1\}$, there is a surjection $f : G \to H$ onto some finite group $H$ such that $f(g) \neq 1$. If the group is countable, this is equivalent to the existence of a descending chain

$$\ldots G_i \subset G_{i-1} \subset \cdots \subset G_0 = G$$

of normal subgroups having finite index such that $\cap_i G_i = 1$.

**Theorem 5:** *Let $X$ be a 4–manifold with even intersection form and residually finite fundamental group. Then the inequality*

$$\frac{5}{4}|\tau(X)| \leq b_2^{(2)}(X)$$

*holds, where $b_2^{(2)}(X)$ denotes the second $L^2$–Betti number of the universal covering.*

*Proof.* We can assume that $X$ is spin, since signature and the $L^2$–Betti numbers are multiplicative in finite coverings. Choose a descending chain of normal subgroups $\Gamma_i \subset \pi_1(X)$ such that $[\pi_1(X) : \Gamma_i]$ is finite and $\cap \Gamma_i = 1$. This sequence corresponds to a tower of coverings $X_i \to X$. Applying Furuta's Theorem to these covering spaces leads to

$$\frac{5}{4}|\tau(X)| \leq \frac{b_2(X_i)}{[\pi_1(X) : \Gamma_i]}.$$

By a Theorem of Lück ([Lu1]), the right hand side of this inequality tends to $b_2^{(2)}(X)$ as $i \to \infty$, and this proves the assertion. □

Note that, since the $L^2$–Euler characteristic equals the usual Euler characteristic, the difference between the second $L^2$–Betti number and $b_2(X)$ depends only on the fundamental group. In particular we obtain the $\frac{5}{4}$–conjecture if

$$b_1^{(2)}(\pi_1(X)) \leq b_1(\pi_1(X)) - 1$$

and $\pi_1(X)$ is residually finite. Therefore Theorem 5 is most interesting if the first $L^2$–Betti number is small or even zero. Some cases where the first $L^2$–Betti number of a group $\Pi$ vanishes are

1. If $\Pi$ contains a normal infinite abelian subgroup (see [CGr]), or

2. if there is an exact sequence

$$1 \longrightarrow \Delta \longrightarrow \Pi \longrightarrow \Gamma \longrightarrow 1$$

with finitely generated groups $\Delta$ and $\Gamma$ such that $\Delta$ is infinite and $\Gamma$ contains $\mathbb{Z}$ as a subgroup (see [Lu2]), or
3. If $\Pi$ is the fundamental group of a manifold $M$ of dimension at least three which is closed and Kähler hyperbolic (see [Gr2]).

In all cases we obtain the inequality $\frac{5}{4}|\tau(X)| \leq e(X)$ (provided that $\pi_1(X) = \Pi$ is residually finite).

Other important cases where the first $L^2$–Betti numbers are completely known are the fundamental groups of compact 3–manifolds which fulfill the Thurston geometrization conjecture [LLu]). These groups are always residually finite (see [He], p. 380).

**Corollary 3:** *Suppose $X$ is a closed 4–manifold with even intersection form whose fundamental group is the fundamental group of a closed hyperbolic manifold of dimension at least three. Then*

$$\frac{5}{4}|\tau(X)| \leq e(X).$$

*In particular, the $\frac{5}{4}$–conjecture holds for $X$ as soon as $b_1(X) > 0$.*

*Proof.* By assumption there is a closed hyperbolic manifold $M$ having fundamental group $\pi_1(X)$. This implies that $\pi_1(X)$ is not finite and is a discrete subgroup of $Iso(\mathbb{H}^n)$, in particular it has a faithful representation into $\mathrm{GL}(n+1;\mathbb{R})$. By a result of Mal'cev (see for example [Lu1]), it is residually finite. Morevover $b_1^{(2)}(X) = b_1^{(2)}(M) = 0$ by a Theorem of Dodziuk ([Dod]), and therefore $e(X) = b_2^{(2)}(X)$. Now the estimate follows from Theorem 5. $\square$

**Remark 4:** Note that the fundamental group of a closed Riemannian manifold with a hyperbolic metric is not amenable. Therefore the last Corollary is in a certain sense contrary to the situation in Theorem 4.

**Remark 5:** The assumption on the dimension in Corollary 3 is really necessary. As an example let $F$ denote a surface of genus 2 with a hyperbolic metric and consider the 4–manifold $X = S^2 \times F$. Then $\pi_1(X)$ is residually finite and the fundamental group of the hyperbolic manifold $F$, but $e(X) = -4$ and $\tau(X) = 0$, hence the inequality $\frac{5}{4}|\tau(X)| \leq e(X)$ does not hold.

Even if the fundamental group of $X$ does not belong to the classes covered by the above results, it is sometimes possible to use similar arguments if the group is obtained as an extension of such a "nice" group by a finitely presentable group. As an example, we mention the following

**Corollary 4:** *Let $X$ be a 4–manifold with even intersection form and assume there is an exact sequence*

$$1 \longrightarrow H \longrightarrow \pi_1(X) \longrightarrow \Gamma \longrightarrow 1$$

*where $H$ and $\Gamma$ are finitely presentable, $\Gamma$ is amenable and $b_1(\Gamma) > 0$. Then*

$$\frac{5}{4}|\tau(X)| \leq e(X).$$

*In particular the $\frac{5}{4}$–conjecture holds for $X$.*

*Proof.* First let us prove the case in which $X$ is spin. Since $b_1(\Gamma) > 0$, we can, using Gromov's argument from [Gr1] as in the proof of Theorem 4, find a sequence $(\Gamma_n)_n$ of normal subgroups $\Gamma_n \subset \Gamma$ with finite index such that $[\Gamma : \Gamma_n] \to \infty$ but $b_1(\Gamma_n)$ is bounded by a constant $C > 0$. Let $G_n \subset \pi_1(X)$ denote the preimages under $\pi_1(X) \to \Gamma$. We obtain exact sequences

$$1 \longrightarrow H \longrightarrow G_n \longrightarrow \Gamma_n \longrightarrow 1$$

and therefore $b_1(G_n) \leq b_1(H) + b_1(\Gamma_n) \leq b_1(H) + C$. Now let $X_n \to X$ denote the covering given by $G_n$, with degree $m_n = [\pi_1(X) : G_n] = [\Gamma : \Gamma_n]$. Using the multiplicativity of the Euler characteristic, we find that

$$\frac{b_2(X_n)}{m_n} = e(X) + \frac{2}{m_n}(b_1(G_n) - 1)$$

$$\leq e(X) + \frac{2}{m_n}(C - 1 + b_1(H)).$$

If we apply Furuta's Theorem to $X_n$, we therefore obtain $\frac{5}{4}|\tau(X)| \leq e(X) + \frac{2}{m_n}(C - 1 + b_1(H))$. The second summand tends to zero if we let $n$ go to infinity and this proves our assertion.

Now let us turn to the general case. By Theorem 1, we can find a finite cyclic cover $\tilde{X} \to X$ such that $\tilde{X}$ is spin. By the multiplicativity of Euler characteristic and signature, we only have to check that the fundamental group of $\tilde{X}$ again satifies our conditions, then we can apply the spin case we just proved to $\tilde{X}$ and we are done. For this purpose, let $G \subset \pi_1(X)$ denote the normal subgroup of finite index corresponding to $\tilde{X}$. Let $\Gamma' \subset \Gamma$ denote its image in $\Gamma$. Then $\Gamma'$ is a subgroup of finite index, hence it is again finitely generated, amenable and has positive first Betti number. If we let $H' = H \cap G$, then $H'$ has finite index in $H$, in particular it is again finitely generated, and we obtain an exact sequence

$$1 \longrightarrow H' \longrightarrow G = \pi_1(\tilde{X}) \longrightarrow \Gamma' \longrightarrow 1$$

as desired. Finally the inequality $\frac{5}{4}|\tau(X)| \leq e(X) = 2 + b_2(X) - 2b_1(X)$ implies the $\frac{5}{4}$–conjecture in this case since $b_1(X) \geq b_1(\Gamma) \geq 1$ by assumption. □

**Example 1:** The condition of Corollary 4 is fulfilled if $b_1(X) > 0$ and the commutator subgroup of $\pi_1(X)$ is finitely presentable, simply let $\Gamma = H_1(X)$ and $H = [\pi_1(X), \pi_1(X)]$.

**Example 2:** Assume that $X$ is a 4–manifold with even intersection form and $\pi_1(X)$ is the fundamental group of a closed manifold $M$ which fibres over the circle. Let $F$ denote the fibre of $M \to S^1$. Then, since $\pi_2(S^1) = 0$, we have an exact sequence
$$1 \longrightarrow \pi_1(F) \longrightarrow \pi_1(M) = \pi_1(X) \longrightarrow \mathbb{Z} \longrightarrow 0$$
and since the fibre is compact, $\pi_1(F)$ is finitely presentable. Hence, by Corollary 4, we obtain $\frac{5}{4}|\tau(X)| \leq e(X)$ and the $\frac{5}{4}$–conjecture for $X$.

## 4. A NEW GEOMETRIC GROUP INVARIANT

In this section, we use the difference $\frac{5}{4}|\tau(\cdot)| - e(\cdot)$, which is an invariant of 4–manifolds, to obtain an invariant for finitely presentable groups (Definition 3), similar to the invariants $q$ (Hausmann–Weinberger, see [HW]) and $p$ (Kotschick, see [Ko1]). We prove some fundamental properties of this invariant and give examples for its computation.

**Definition 3:** *Let $\Gamma$ be a finitely presentable group. We define*
$$r(\Gamma) = inf\{e(X) - \frac{5}{4}|\tau(X)| \mid X \text{ even }, \pi_1(X) = \Gamma\}.$$
*where the infimum is taken over all closed, connected and oriented smooth 4–manifolds with even intersection form and fundamental group $\Gamma$.*

**Lemma 7:** *For every finitely presentable group $\Gamma$, $r(\Gamma)$ is an even integer with*
$$(1) \qquad r(\Gamma) \geq 2(1 - dim_{\mathbb{Z}_2}H^1(\Gamma;\mathbb{Z}_2)).$$
*In particular, there is an even 4–manifold $X$ with $e(X) - \frac{5}{4}|\tau(X)| = r(\Gamma)$ and $\pi_1(X) = \Gamma$, i.e. the infimum is a minimum. Equality in (1) occurs if and only if $\Gamma$ is the fundamental group of a 4–manifold $X$ with $H^2(X;\mathbb{Z}_2) = 0$.*

*Proof.* Suppose that $X$ is a 4–manifold with even intersection form and fundamental group $\Gamma$. We have
$$e(X) - \frac{5}{4}|\tau(X)| = dim_{\mathbb{Z}_2}H^2(X;\mathbb{Z}_2) + 2 - 2\,dim_{\mathbb{Z}_2}H^1(X;\mathbb{Z}_2) - \frac{5}{4}|\tau(X)|$$
$$= 2 - (\frac{5}{4}|\tau(X)| - dim_{\mathbb{Z}_2}H^2(X;\mathbb{Z}_2)) - 2\,dim_{Z_2}H^1(X;\mathbb{Z}_2)$$
and we know from Theorem 3 that $\frac{5}{4}|\tau(X)| - dim_{\mathbb{Z}_2}H^2(X;\mathbb{Z}_2) \leq 0$. This implies the claimed inequality. In particular the infimum is not $-\infty$ and must be a minimum. Now assume that $r(\Gamma) = 2(1 - dim_{\mathbb{Z}_2}H^1(\Gamma;\mathbb{Z}_2))$. Choose $X$ realizing the minimum. Then the above calculation shows that $\frac{5}{4}|\tau(X)| - dim_{\mathbb{Z}_2}H^2(X;\mathbb{Z}_2) = 0$. By Theorem 3, this implies $\tau(X) = 0$ and hence $H^2(X;\mathbb{Z}_2) = 0$. Conversely if there is a 4–manifold $X$ with $H^2(X;\mathbb{Z}_2) = 0$ and $\pi_1(X) = \Gamma$, we can deduce $b_2(X) = 0$, in particular the intersection form of $X$ is even. This shows that $r(\Gamma) \leq e(X) = 2(1 - dim_{\mathbb{Z}_2}H^1(X;\mathbb{Z}_2)) \leq r(\Gamma)$, and we obtain equality. □

**Remark 6:** If we let the infimum go over all $X$ with $\pi_1(X) = \Gamma$, not only the ones with even intersection form, we obtain $-\infty$, because one can always add a copy of $\mathbb{C}P^2$ and thereby decrease $e(\cdot) - \frac{5}{4}\tau(\cdot)$ by $\frac{1}{4}$. Also note that one of the standard constructions for a 4–manifold with prescribed $\pi_1$ – realizing this group as the fundamental group of a 2–complex, embedding this into $\mathbb{R}^5$ and taking the boundary $X$ of a regular neighborhood – provides a 4–manifold which is spin (since $TX \oplus \mathbb{R}$ is the tangent bundle of $\mathbb{R}^5$ and hence trivial), therefore $r(\Gamma)$ is defined for all finitely presentable groups $\Gamma$. Another standard argument that every finitely presented group is the fundamental group of a spin 4–manifold is the following: realise the generators by gluing copies of $S^1 \times S^3$ and then perform surgery along circles representing the relations, one easily checks that the surgery can be done in a spin–preserving way.

**Lemma 8:** *If $\Gamma' \subset \Gamma$ is a subgroup of finite index $[\Gamma, \Gamma']$, then*

$$r(\Gamma') \leq [\Gamma, \Gamma']\, r(\Gamma).$$

*Furthermore, we have the estimate*

$$r(\Gamma_1 * \Gamma_2) \leq r(\Gamma_1) + r(\Gamma_2) - 2.$$

*Proof.* The first equation is a consequence of the multiplicativity of $\tau$ and $e$ under finite covers: Choose a 4–manifold $X$ with even intersection form and $\pi_1(X) = \Gamma$ such that $e(X) - \frac{5}{4}|\tau(X)| = r(\Gamma)$. Then the subgroup $\Gamma'$ is the fundamental group of a finite cover $\tilde{X} \to X$ which again has even intersection form. Now $e(\tilde{X}) = [\Gamma, \Gamma']\, e(X)$, $\tau(\tilde{X}) = [\Gamma, \Gamma']\, \tau(X)$, and therefore $r(\Gamma') \leq [\Gamma, \Gamma']\, r(\Gamma)$ as claimed.

Now suppose $X_1$ and $X_2$ are even 4–manifolds with fundamental groups $\Gamma_1$ and $\Gamma_2$ such that $e(X_i) - \frac{5}{4}|\tau(X_i)| = r(X_i)$. Orient both manifolds such that their signatures are non–negative. Then $X_1 \# X_2$ has even intersection form, fundamental group $\Gamma_1 * \Gamma_2$, signature $\tau(X_1) + \tau(X_2)$ and Euler characteristic $e(X_1) + e(X_2) - 2$, and the assertion follows. $\square$

In some cases, the invariant $r$ can be computed or at least estimated using the results of the last section. We summarize some of these computations in the following proposition.

**Proposition 1:** *Let $\Gamma$ be a finitely presentable group.*

1. *If $\Gamma$ is free of rank $k$, then $r(\Gamma) = 2 - 2k$.*
2. *If $\Gamma$ is a cyclic group of finite order, then $r(\Gamma) = 2$.*
3. *For every finitely presentable group $\Gamma$, we have $r(\Gamma * \mathbb{Z}) = r(\Gamma) - 2$.*
4. *If $\Gamma$ is amenable, then $r(\Gamma) \geq 0$.*
5. *If $\Gamma$ is the fundamental group of a closed oriented surface of positive genus $g$, then $r(\Gamma) = 4 - 4g$.*
6. *We have $r(\mathbb{Z}^2) = r(\mathbb{Z}^4) = 0$ and $0 \leq r(\mathbb{Z}^3) \leq 2$.*

*Proof.* 1. Let $F_k$ denote the free group of rank $k$. Consider the manifolds $X_0 = S^4$ and $X_k = k(S^1 \times S^3)$. Then $\pi_1(X_k)$ is free of rank $k$, $X_k$ has even intersection form (actually it is spin) and $e(X_k) = 2 - 2k, \tau(X_k) = 0$. This shows $r(F_k) \leq 2 - 2k$. But from Lemma 7 we know $r(F_k) \geq 2 - 2k$, hence the claimed equality holds.

2. Let $n = \#\Gamma$. Then $\Gamma = \mathbb{Z}_n$, and there is a rational homology 4–sphere $X$ with fundamental group $\Gamma$, which can be constructed by starting with $S^1 \times S^3$ and performing surgery along a circle representing $n$ times the generator of $\pi_1(S^1 \times S^3)$. We have $\tau(X) = 0$ and $e(X) = 2$, therefore $r(\Gamma) \leq 2$. On the other hand, Corollary 2 implies that for every 4–manifold $Y$ with even intersection form and fundamental group $\Gamma$, $\frac{5}{4}|\tau(Y)| \leq b_2(Y)$ holds, and $b_2 = e(Y) - 2$, hence $e(Y) - \frac{5}{4}|\tau(Y)| \geq 2$.

3. We can use exactly the same arguments as in [Ko2], Theorem 3.6: Choose an even 4–manifold $X$ with $\pi_1(X) = \Gamma * \mathbb{Z}$ and $e(X) - \frac{5}{4}|\tau(X)| = r(\Gamma * \mathbb{Z})$. We can do surgery along a generator of $\mathbb{Z} \subset \pi_1(X)$ to obtain a 4–manifold $Y$ with fundamental group $\Gamma$ which has the same intersection form as $X$. Then $e(Y) = e(X) + 2$, and therefore $e(Y) - \frac{5}{4}|\tau(Y)| = r(\Gamma * \mathbb{Z}) + 2$. This proves $r(\Gamma) \leq r(\Gamma * \mathbb{Z}) + 2$. Now Lemma 8 and $r(\mathbb{Z}) = 0$ imply the claimed equality.

4. This is an immediate consequence of Theorem 4.

5. The proof is essentially the same as the proof for the corresponding property of the invariant $p$ in [Ko2]. We repeat the argument here for sake of completeness.

Let $\Gamma_g$ denote the fundamental group of the closed oriented surface of genus $g$. We have $b^1(\Gamma_g; \mathbb{Z}_2) = 2g$, and by Lemma 7, this implies $r(\Gamma_g) \geq 2 - 4g$. Furthermore the example of the spin–manifold $F_g \times S^2$ shows $r(\Gamma_g) \leq 4 - 4g$. Since the invariant is even, we only have to rule out the case $r(\Gamma_g) = 2 - 4g$.

Suppose that there is in fact some $g \geq 1$ with $r(\Gamma_g) = 2 - 4g$. Choose a spin 4–manifold $X$ with $\pi_1(X) = \Gamma_g$ and $e(X) - \frac{5}{4}|\tau(X)| = 2 - 4g$. Consider a double cover $F_h \to F_g$. By multiplicativity of the Euler characteristic, we have $h = 2g - 1$. The same cover will define a double cover $X' \to X$, and $e(X') = 2e(X), \tau(X') = 2\tau(X)$. Therefore we obtain the inequality

$$e(X') - \frac{5}{4}|\tau(X')| = 2(2 - 4g) \geq r(\Gamma_h) \geq 2 - 4h$$

and hence $4 - 8g \geq 2 - 4h = 6 - 8g$, which is a contradiction.

6. In all the three cases, the group is amenable and hence $r(\Gamma) \geq 0$. If $\Gamma = \mathbb{Z}^2$ or $\Gamma = \mathbb{Z}^4$, the examples $T^2 \times S^2$ and $T^2 \times T^2$ show that $r \leq 0$. Furthermore we can surger a circle on $T^4$ to obtain a spin 4–manifold with signature zero, Euler characteristic 2 and fundamental group $\mathbb{Z}^3$. Hence $r(\mathbb{Z}^3) \leq 2$. □

The following table contains all groups for which the three invariants $p, q$ and $r$ are known (at least to the author). The values for $r$ are taken from the preceeding proposition, for the claims on $p$ and $q$, we refer the reader to [Ko2]. In the table, $F_r$ denotes the free group of rank $r$, and $\Gamma_g$ denotes the fundamental group of a closed orientable surface of genus $g$.

| $\Gamma$ | $p$ | $q$ | $r$ |
|---|---|---|---|
| 1 | 2 | 2 | 2 |
| $F_r$ | $2-2r$ | $2-2r$ | $2-2r$ |
| $\mathbb{Z}_n$ | 2 | 2 | 2 |
| $\Gamma' * \mathbb{Z}$ | $p(\Gamma')-2$ | $q(\Gamma')-2$ | $r(\Gamma')-2$ |
| $\mathbb{Z}^2$ | 0 | 0 | 0 |
| $\mathbb{Z}^4$ | 0 | 0 | 0 |
| $\mathbb{Z}^3$ | 2 | 2 | ? |
| $\Gamma_g$ | $4-4g$ | $4-4g$ | $4-4g$ |

**Question:** What is $r(\mathbb{Z}^3)$?

The table above might suggest that the three invariants $p, q$ and $r$ agree. However, this is not the case. In [Ko2], D. Kotschick proved that the invariants $p$ and $q$ are in general different. We will now see that there is in fact an infinite family of groups for which all the three invariants are different.

**Theorem 6:** *There is a sequence $(\Gamma_n)_{n\geq 1}$ of finitely presentable groups such that $b_1(\Gamma_n) \to \infty$ and*

$$r(\Gamma_n) < p(\Gamma_n) < q(\Gamma_n).$$

*Proof.* We use the examples $M_{n,m} \to S$ for aspherical surface bundles over surfaces with positive signatures $\tau(M_{n,m})$ constructed by Kodaira in [Kod]. The manifold $M_{n,m}$ is an algebraic surface which is an m–fold cover over a product of algebraic curves, branched along a divisor. The formula for the canonical class of $M_{n,m}$ in [Kod] immediately shows that $M_{n,m}$ is spin if $m$ is odd. Hence we can construct a family $X_n$ of finite covers of $M_{2,3}$ with $b_1(X_n) \to \infty$ by taking the pullbacks of $M_{2,3} \to S$ to finite covers of $S$ (note that $g(S) > 1$). Clearly each $X_n$ is aspherical, spin and has positive signature. Furthermore, $e(M_{2,3}) > 0$ and therefore $e(X_n) > 0$. Let $\Gamma_n = \pi_1(X_n)$. By Theorem 3.8 in [Ko2], $X_n$ realizes $p$ and $q$, i.e. $e(X_n) - |\tau(X_n)| = p(\Gamma_n) < q(\Gamma_n) = e(X_n)$. Moreover we have

$$r(\Gamma_n) \leq e(X_n) - \frac{5}{4}|\tau(X_n)| < p(\Gamma_n),$$

since $X_n$ is spin. This implies $r(\Gamma_n) < p(\Gamma_n) < q(\Gamma_n)$, as desired. $\square$

Department of Mathematics, Yale University, P.O. Box 208283, New Haven, CT 06520 – 8283, USA

*E-mail address*: bohr@math.yale.edu